\documentclass[11pt]{article}
\usepackage{times,amssymb,amsfonts}
\def\qed{\ifhmode\unskip\nobreak\fi\ifmmode\ifinner\else\hskip.5em\fi\fi
 \hbox{\hskip.5em$\square$\hskip.1em}}
\newenvironment{genericem}[1]{\smallskip{\sc #1.}\em}{\/\rm \smallskip}
\newenvironment{proof}{\smallskip{\sc Proof.}}{\qed\smallskip}
\begin{document}
\title{Stable Configurations of repelling Points on compact hyperbolic Manifolds}
\author{Burton Randol}
\date{}
\maketitle

In this note, we describe a feature of stable configurations of repelling points on certain manifolds in the sense of \cite{randolchapter}\ (pp.\ 282--288). The approach is applicable in wider contexts than are discussed here, but we will confine our discussion to compact hyperbolic manifolds, since in this case the necessary preparatory material is already available in \cite{randolchapter}. The discussion of this subject in \cite{randolchapter} was motivated by analogies with classical electrostatics, but the conclusions drawn there can, if preferred, be regarded as being of a formal character, and independent of these analogies.

Accordingly, supposing $M$ to be a compact hyperbolic manifold, let $S$ be a finite set of distinct points in $M$. We regard the points of $S$ as mutually repelling along connecting geodesics, with the repelling magnitude varying with distance, as specified by a function $H(\rho)$ $\rho \in [0,\infty)$. In other words, the effect exerted on a point $p \in S$ by points in $S$ is computed by taking the (infinite) vector sum $\sum_n H(L_n)V_n$, where $L_n$ ranges over the lengths of the geodesic segments $g_n$ connecting $p$ to points of $S$, and $V_n$ is the unit tangent vector to $g_n$ at $p$, taken in the repelling direction. Note that $p$ is permitted to act on itself. Conceptually, this can be visualized intrinsically on $M$, or if one prefers, on the universal covering space $H^n$, where the mutual repulsions act along the countably numerous geodesics connecting points, attenuated according to $H(\rho)$. Similar considerations can also be applied to continuous, rather than point distributions on $M$, where the repulsion between infinitesimal parts of $M$ is defined analogously \cite{randolchapter}, but in this note we will be concerned with finite point distributions. In the interests of expository economy and to avoid unnecessary duplication, we refer the reader to \cite{randolchapter} for background material concerning the Selberg pretrace and trace formulas, as well as to the application of the pretrace formula to the matters considered here. In greater detail, following the discussion in \cite{randolchapter}, we assume that $H(\rho)$ is linked to a function $k(\rho)$ by the requirement that $H(\rho) = -k'(\rho)$, where $k(\rho)$, which corresponds to a point-pair invariant, is a smooth function on $[0,\infty)$ which vanishes at infinity, and for which the Selberg pretrace formula is valid, with uniform and absolute convergence. As is customary, we denote the Selberg transform of $k(\rho)$ by $h(r)$, and we will be concerned with the values of $h(r)$ on the cross-shaped subset of the complex plane composed of the union $L$ of the real axis with the closed segment of the imaginary axis from \mbox{$-(n-1)i/2$} to $(n-1)i/2$. In the standard parameterization $r \leftrightarrow \lambda$, this corresponds to eigenvalues of the positive Laplacian on $M$, which are situated on $[0,\infty)$. We will also impose, for the purposes of the present note, the additional requirement that $h(r)$ is positive on $L$, and will call the configuration $S$ stable, if the net repelling effect at each point of $S$ is zero.

Let $0=\lambda_0 < \lambda_1 \leq \lambda_2 \ldots$ be the eigenvalues of the positive Laplacian on $M$, with repetitions to account for multiplicity, and suppose $\varphi_0, \varphi_1, \varphi_2,\ldots$ is an associated orthonormal sequence of eigenfunctions. Denote the points of the set $S$ by $x_1,\ldots , x_N$. With this notation, it follows from the discussion in \cite{randolchapter} (cf.\ especially p. 287), that sets $S$ for which the quantity \[\sum_{n=1}^{\infty} \left(\sum_{i=1}^N \sum_{j=1}^N \varphi_n(x_i)\overline{\varphi_n}(x_j)\right)h(r_n)\] \[=  \sum_{n=1}^{\infty} |\varphi_n(x_1) + \cdots \varphi_n(x_N)|^2 h(r_n)\] is locally minimized in the $N$-fold Cartesian product of $M$ with itself are stable in the above sense. In somewhat more detail, the gradient of the last expression, which must vanish at a minimum, can be analyzed using the formula (17) from page $283$ of \cite{randolchapter}, which describes the repelling effect produced by a point of $M$ on another point of $M$. 


This suggests, as being of particular interest, the study of the statistical behavior, as $N \rightarrow \infty$, of point configurations that are globally minimizing. and as we will now show, such configurations are equidistributed in the limit. In order to see this, we will prove a somewhat more general fact, expressed as a theorem, from which the assertion immediately follows.

\begin{genericem}{Theorem} Suppose that $M$ is a compact Riemannian manifold,
and that $1,\varphi_1,\varphi_2,\ldots$ is a sequence of
continuous functions in $C(M)$, each having $L^2$ norm $1$, and with $\int_M \varphi_i = 0$ for $i=1,2,\ldots$. Suppose also that finite linear combinations of functions from the sequence are sup-norm dense in $C(M)$. Denote by 
$M^N$  the $N$-fold Cartesian product of $M$ with itself. Let
$a_1,a_2,\ldots$ be a sequence of positive numbers for which the series
$\sum_{n=1}^\infty |\varphi_n(x)|^2 \, a_n$ is uniformly convergent, and
suppose, for $N=1,2,3\ldots$, that $X_N \in M^N$ globally minimizes the
quadratic expression \[ \sum_{n=1}^\infty \left| \varphi_n(x_1)+ \ldots
+\varphi_n(x_N)\right| ^2 \, a_n \,.\] Then, as $N \rightarrow \infty$,
the sets $X_N$ become equidistributed, in the sense that for any
continuous function $f$ on $M$, \[\frac{1}{N}\sum_{x \in X_N} f(x)
\rightarrow \int_M f\] as $N \rightarrow \infty$.\end{genericem}

\begin{proof} It follows from the hypotheses that for each $n$, \[
\int_{M\times\cdots\times M} \left| \varphi_n(x_1)+ \ldots +
\varphi_n(x_N)\right| ^2 \,dx_1\cdots dx_k = NV^{N-1} \,,\]where $V$ is the volume of $M$, so
\[\int_{M\times\cdots\times M} \left(\sum_{n=1}^\infty \left|
\varphi_n(x_1)+ \ldots + \varphi_n(x_N)\right| ^2 a_n\right)
\,dx_1\cdots dx_k = NV^{N-1}\sum_{n=1}^\infty a_n\,. \]

In particular, by the mean value theorem for integrals, there exists $(\overline{x}_1,\ldots ,\overline{x}_N)\in
M^N$ such that \[ \sum_{n=1}^\infty \left| \varphi_n(\overline{x}_1)+
\ldots + \varphi_n(\overline{x}_N)\right| ^2 a_n = \frac{N}{V} \sum_{n=1}^{\infty}
a_n\,,\] which shows that if $(p_1\ldots ,p_N) \in M^N$ is globally
minimizing, then \[ \sum_{n=1}^\infty \left| \varphi_n(p_1)+ \ldots +
\varphi_n(p_N)\right| ^2 a_n \leq \frac{N}{V} \sum_{n=1}^{\infty} a_n \,.\]

For a fixed $m$, we therefore conclude from the positivity of the $a_n$'s that \[ \left| \varphi_m(p_1)+
\ldots + \varphi_m(p_N)\right|^2 \leq \frac{N}{V} a_m^{-1} \sum_{n=1}^{\infty}
a_n \,,\] so \[ \left| \varphi_m(p_1)+ \ldots + \varphi_m(p_N)\right|
\leq \sqrt{\frac{N}{V} a_m^{-1} \sum_{n=1}^{\infty} a_n} \;\;,\] which implies
that \[ N^{-1} \left| \varphi_m(p_1)+ \ldots + \varphi_m(p_N)\right|
\leq \frac{C(m)}{\sqrt{N}}\,,\] where \[ C(m) = \sqrt{\frac{1}{Va_m}
\sum_{n=1}^{\infty} a_n} \;\;.\] In particular, the $X_k$'s in the limit integrate 
the uniformly dense set of functions $1,\varphi_1,\varphi_2,\ldots$ (the constant function is automatic), so Weyl's criterion is satisfied, and the $X_k$'s are equidistributed in the limit.\end{proof}

\bigskip

\noindent\textbf{Remarks.}\begin{enumerate}

\item If we consider globally minimizing configurations to be equivalent when they are mutual images by isometries of $M$, it would be interesting to know something about the behavior with $N$ of the number of inequivalent configurations.

\item Our definition of point interaction is intrinsic, and does not involve an isometric embedding of $M$ in, e.g., a Euclidean space. A prominent example of the latter occurs in the well-known Thompson problem concerning stable configurations of electrons on the unit sphere of $R^3$, where the repelling actions are transmitted via mutually connecting geodesic segments in the containing space $R^3$. In view of various isometric embedding theorems, it would be interesting to explore connections between the intrinsic and extrinsic viewpoints.

\end{enumerate}

\newcommand{\noopsort}[1]{}

\bigskip

\begin{flushleft}
{\sc Ph.D Program in Mathematics\\CUNY Graduate Center\\365 Fifth
Avenue\\New York, NY 10016}
\end{flushleft}

\end{document}